\documentclass[reqno]{amsart}

\usepackage[english]{babel}
\usepackage[utf8x]{inputenc}
\usepackage{textcase}

\usepackage[pdftex,pdfcreator={LaTeX},pdfproducer={PDFLaTeX},colorlinks,
	bookmarks,
	pdftitle={},
	pdfauthor={},
	pdfkeywords={},
	linkcolor=blue,                  % color of internal links
	citecolor=red,                   % color of links to bibliography
	filecolor=red,                   % color of file links
	urlcolor=cyan                    % color of external links
]{hyperref}

\usepackage{environ}
\usepackage{enumerate}

\NewEnviron{highlights}{}
\NewEnviron{keyword}{}

\usepackage{amsmath,amsthm,amsfonts,mathtools}
\usepackage{longtable}
\usepackage[colorlinks]{hyperref}
\usepackage{kbordermatrix}

\newcommand{\G}{\ensuremath{\Gamma}}
\newcommand{\Z}{\ensuremath{\mathbb{Z}}}

\newcommand{\dbab}{\ensuremath{(d, b, \alpha, \beta)}}

\newcommand{\doi}[1]{\href{https://doi.org/#1}{\texttt{doi:#1}}}

\newcommand{\noop}[1]{}

\newtheorem{theorem}{Theorem}[section] % 1st argument is your name for it
\newtheorem{lemma}[theorem]{Lemma}     % 2nd argument is what is printed
\newtheorem{corollary}[theorem]{Corollary}
\newtheorem{claim}[theorem]{Claim}
\newtheorem{proposition}[theorem]{Proposition}
\newtheorem{remark}[theorem]{Remark}

\newcommand{\arxivbr}{\\}

\begin{document}
\gdef\papertitle{On local spectra of classical distance-regular graphs with $b = 2$}

\gdef\authoronename{Alexander L.~Gavrilyuk}
\newcommand{\authoroneabbrv}{A.L.\ Gavrilyuk}
\newcommand{\authoroneorg}{University of Memphis}
\newcommand{\authoronestate}{Tennessee}
\newcommand{\authoronecountry}{United States of America}
\newcommand{\authoroneemail}{A.Gavrilyuk@memphis.edu}

\gdef\authortwoname{Jano\v{s} Vidali}
\newcommand{\authortwoabbrv}{J.\ Vidali}
\newcommand{\authortwoorgone}{Faculty of Mathematics and Physics, University of Ljubljana}
\newcommand{\authortwoorgtwo}{Institute of Mathematics, Physics and Mechanics}
\newcommand{\authortwocity}{Ljubljana}
\newcommand{\authortwocountry}{Slovenia}
\newcommand{\authortwoemail}{janos.vidali@fmf.uni-lj.si}

\newcommand{\dedication}{%
Dedicated to Professor Jack Koolen on the occasion of his 60th birthday \arxivbr
and to Professor Akihiro Munemasa on the occasion of his retirement}

\title{\papertitle}

\author[\authoroneabbrv]{\authoronename}
\address{\authoronename,\newline
	\authoroneorg, \authoronestate, \authoronecountry.
}
\email{\authoroneemail}

\author[\authortwoabbrv]{\authortwoname}
\address{\authortwoname,\newline
	\authortwoorgone, \authortwocountry.
	\newline also affiliated with: \newline
	\authortwoorgtwo, \authortwocity, \authortwocountry.
}
\email{\authortwoemail}

\dedicatory{\dedication}

%% Abstract
\begin{abstract}
%% Text of abstract
Using a combination of spectral techniques and combinatorial counts,
we prove a number of non-existence results
for putative distance-regular graphs with classical parameters,
focusing on the case $b=2$.
Among those, the case $(d,b,\alpha,\beta)=(d,2,1,2^d-2)$
passed all known feasibility conditions
and required a deeper combinatorial analysis,
from which nonexistence has been derived whenever $d \ge 9$.
In doing so,
we constructed a new family of regular graphs with four distinct eigenvalues
which would have been local graphs in such a distance-regular graph.
\end{abstract}

%%Graphical abstract
%\begin{graphicalabstract}
%\includegraphics{grabs}
%\end{graphicalabstract}

%%Research highlights
%\begin{highlights}
%\item Research highlight 1
%\item Research highlight 2
%\end{highlights}

%% Keywords
\begin{keyword}
%% keywords here, in the form: keyword \sep keyword

distance-regular graphs \sep Terwilliger polynomial \sep classical parameters

%% PACS codes here, in the form: \PACS code \sep code

%% MSC codes here, in the form: \MSC code \sep code
%% or \MSC[2008] code \sep code (2000 is the default)

\MSC[2020] 05E30

\end{keyword}

\maketitle

\section{Introduction}

One of the ultimate problems in algebraic combinatorics is to classify the families of distance-regular graphs with unbounded diameter~\cite{BanIto1984}.
All such \emph{known} families\footnote{%
Except the polygons, i.e., the cycles.
}
either have the so-called classical parameters
or can be derived from the former ones~\cite[Chapter~6]{BCN}, \cite[Section~3.1]{vanDamKoolenTanaka2016},
which thus will be referred to as classical graphs.
We will give the exact definitions in Section~\ref{sec:clas};
for now we just say that the classical parameters
associated with a classical graph
are written as a $4$-tuple $(d,b,\alpha,\beta)$;
here, $d$ is simply the diameter of the graph in question,
and $b$ is an integer not equal to $0$ or $-1$.
Apart from the well-known families,
there are many sets of classical parameters
for which the existence of a corresponding graph is still open.

Terwilliger~\cite[Theorem~6.1.1]{BCN} classified
classical distance-regular graphs having $b=1$,
which include well-known examples such as Hamming graphs and Johnson graphs. Weng~\cite{Weng1999} (see also~\cite{Tian2024})
classified classical graphs with $b<-1$ under some mild conditions.

Two other prominent families of classical distance-regular graphs,
the bilinear forms graphs and Grassmann graphs,
have $b\geq 2$ and they turn out to be the only examples\footnote{%
The parameter $b$ is equal to the order $q$ of the finite field
over which the graph in question is defined.
}
if the parameter $\beta$ is sufficiently large
in terms of the other parameters $\alpha$, $b$, and $d$,
as shown by Metsch~\cite{Metsch1999}.
Thus,
a lot of research~\cite{JV17,KOOLEN2026104378,KoolenLvParkYang2025,LvKoolen2026}
has been devoted to improving a lower bound for $\beta$
and finding other bounds and feasibility conditions on the classical parameters.

The focus of this paper is on the case $b=2$,
as a natural next step after $b=1$.
It is known that the smallest eigenvalue of a \emph{local graph},
i.e., the subgraph induced by the neighborhood of a vertex in a classical graph with parameter $b$,
is at least $-b-1$.
This allowed Terwilliger to invoke a celebrated result~\cite{Cameron1976}
on graphs with smallest eigenvalue at least $-2$ for $b=1$.
Koolen et al.~\cite{KOOLEN2026104378,  KoolenJYYangQYang2018,YangKoolen2025}
developed a structure theory for graphs with smallest eigenvalue at least $-3$,
but this theory is currently too complex to be readily applied
to complete the classification for the $b=2$ case.

We make use of another result of Terwilliger~\cite{GK15,Suzuki}
which states that the eigenvalues of the local graphs
belong to the non-negativity intervals of some degree $4$ polynomial.
The coefficients and hence the roots of this polynomial
are completely determined by the classical parameters $(d,b,\alpha,\beta)$.
For certain values of classical parameters,
these intervals become narrow enough
so that we can determine the spectrum of local graphs
by leveraging the fact that non-integral eigenvalues must be algebraic integers, while their multiplicities are controlled by standard combinatorial counts.
This gives a number of non-existence results, not limited to the case $b=2$,
that we present in this paper.

This spectral technique generalizes the idea from~\cite{GK19}
where it was used to prove that the bilinear $(d\times d)$-forms graphs
over the binary field $\mathbb{F}_2$
are the only graphs with classical parameters $(d,2,1,2^{d}-1)$.
Recall that the bilinear $(d\times e)$-forms graph, $e\geq d$,
is defined on the set of $(d\times e)$-matrices
over the finite field $\mathbb{F}_q$,
with two matrices being adjacent if the rank of their difference is equal to 1;
it has classical parameters $(d,q,q-1,q^{e}-1)$.

We have some evidence suggesting
that the bilinear forms graphs over $\mathbb{F}_2$
might be the only distance-regular graphs
with classical parameters $(d,2,1,\beta)$.
The above-mentioned results from~\cite{Metsch1999} and~\cite{GK19}
confirm this if $\beta\geq 2^{d+4}-1$ and $\beta=2^{d}-1$,
respectively.
On the other hand, note that $\beta\geq 2^{d-1}$ holds
by~\cite[Proposition~1(iii)]{JV17}.

To make use of our spectral technique, we assume that $\beta\leq 2^d-1$,
in which case we get a number of non-existence results for classical parameters $(d,2,1,\beta)$
which follow from standard feasibility conditions. In doing so,
we found the parameters $(d,2,1,2^{d}-2)$
which passed all feasibility conditions.
Assume $\Gamma$ is a classical graph with such parameters. We were able to determine the spectrum of local graphs in $\Gamma$,
which also did not reveal any contradiction.
We went further and constructed a family of graphs
having the prescribed spectrum.
Using the technique developed by Koolen et al.~in~\cite{GK19,grids,thinGrassmann},
we then show that the local graphs in $\Gamma$
must be isomorphic to graphs from our construction
if the diameter $d$ of $\Gamma$ is at least~$9$.
While it did briefly appear that we might have been on the verge of
constructing a new family of distance-regular graphs with unbounded diameter,
unfortunately,
a deeper combinatorial analysis shows that
such graphs do not exist whenever $d \ge 9$,
and we have reasons to believe that this is the case for all $d\geq 3$.
In fact,
the authors are aware of an AI-generated nonexistence proof for $d \ge 3$;
it is based on the fact that $\Gamma$ must be
a so-called 1-thin distance-regular graph,
which follows from the spectrum of its local graphs.
We defer the analysis and possible generalization of this proof to future work.

Furthermore, it was shown by Terwilliger~\cite[Theorem~77]{Suzuki} that a stronger property,
\emph{thinness},
implies that the local graphs in classical graphs with diameter $d\geq 5$
must be \emph{co-edge-regular}~\cite{coedge},
meaning that each pair of non-adjacent vertices
has the same number of common neighbors.
This property does not hold in our case;
thus, another possible way to rule out the existence of graphs
with classical parameters $(d,2,1,2^d-2)$
is to generalize~\cite[Theorem~77]{Suzuki}
under a weaker condition of $1$-thinness.

In conclusion, our (unfortunately, unsuccessful) attempt to construct a new family of graphs with classical parameters $(d,2,1,2^{d}-2)$ underscores the complexity of the complete classification of classical graphs with $b=2$, which seems to be elusive at the moment.

Note that our construction of graphs
which would have been local graphs in $\Gamma$ is of independent interest. 
Graphs on at most $30$ vertices having the prescribed spectrum
were known in the literature~\cite{vD95,vanDamSpence1998} from computer search,
but the first candidate with $d=3$ in the putative family
with classical parameters $(d,2,1,2^{d}-2)$ would have valency $42$.
The corresponding (local) spectrum $\{[11]^1,[4]^{12},[-1]^{14},[-3]^{15}\}$
was known to experts~\cite{vDemail} as feasible,
yet no actual graph constructions were known.

\section{Spectra of local graphs in $Q$-polynomial distance-regular graphs} \label{sec:spectra}

We follow the standard terminology and notation in the theory of distance-regular graphs, see, e.g., \cite{BCN,vanDamKoolenTanaka2016}.
Let $\G$ be a distance-regular graph of diameter $d$
with intersection array $\{k, b_1, \dots, b_{d-1}; 1, c_2, \dots, c_d\}$,
and set $b_0 = k$, $b_d = c_0 = 0$, $c_1 = 1$,
and $a_i = k - b_i - c_i$ ($0 \le i \le d$).
Clearly, for each vertex $u$ of $\G$,
the local graph $\G(u)$ is an $a_1$-regular graph with $k$ vertices.
We will use the following results
to obtain information on the spectra of these local graphs.

\begin{theorem}\label{thm:bounds} {\rm \cite[Theorem~4.4.3(i)]{BCN}}
Let $\G$ be a distance-regular graph of diameter $d \ge 3$
with distinct eigenvalues $k = \theta_0 > \theta_1 > \cdots > \theta_d$,
and let $u$ be a vertex of $\G$.
Denote the eigenvalues of the local graph $\G(u)$
by $a_1 = \eta_1 \ge \eta_2 \ge \cdots \ge \eta_k$.
Then $\eta_2 \le -1 - {b_1 \over \theta_d + 1} =: \eta_{\max}$
and $\eta_k \ge -1 - {b_1 \over \theta_1 + 1} =: \eta_{\min}$.
\qed
\end{theorem}

If the graph $\G$ is $Q$-polynomial,
then the so-called Terwilliger polynomials $T_i$
($2 \le i \le d-1$)~\cite[\S4]{GK15}, \cite{Suzuki}
can be used to further restrict
the non-principal eigenvalues of its local graphs.

\begin{theorem}\label{thm:terwpoly} {\rm \cite[Theorem~4.2]{GK15}}
Let $\G$ be a $Q$-polynomial distance-regular graph of diameter $d \ge 3$,
and let $u$ be a vertex of $\G$.
Then,
for every eigenvalue $\eta$ of the local graph $\G(u)$
such that $\eta \ne k$,
$T_i(\eta) \ge 0$ holds for all $2 \le i \le d-1$.
\qed
\end{theorem}

Since the Terwilliger polynomials $T_i$ ($2 \le i \le d-1$)
are degree $4$ polynomials which are independent of $i$ up to a scalar multiple,
Theorems~\ref{thm:bounds} and~\ref{thm:terwpoly}
can be used to limit the non-principal eigenvalues of the local graphs
of a $Q$-polynomial distance-regular graph of diameter at least $3$
to a union of at most three bounded intervals.
Moreover, in case there exist some indices $i, j$ ($2 \le i, j \le d-1$)
such that $T_i$ and $T_j$ differ in the signs of the leading coefficients,
then the non-principal eigenvalues of the local graphs
are limited to their roots.

In the case when all non-integral eigenvalues of the local graphs
can be shown to lie in a narrow enough interval,
we can further restrict their possible values.

\begin{lemma}\label{lem:noninteig}
Let $\Delta$ be a graph such that all of its nonintegral eigenvalues
lie in an interval $(m - \delta, m + \epsilon)$ for some integer $m$
and some real numbers $\delta, \epsilon$
such that $0 < \delta \le 1$ and $0 < \epsilon \le 2$.
If $\eta$ is a nonintegral eigenvalue of $\Delta$,
then $\delta \ge {-1 + \sqrt{5} \over 2}$, $\epsilon \ge {1 + \sqrt{5} \over 2}$
and $\eta = m + {1 \pm \sqrt{5} \over 2}$.
\end{lemma}

\begin{proof}
Since the eigenvalues of a graph are roots
of the characteristic polynomial of its adjacency matrix,
which is monic with integral coefficients,
it follows that eigenvalues are real algebraic integers,
and for each eigenvalue,
its conjugates appear as eigenvalues with the same multiplicity.

Suppose that $\eta \in (m - \delta, m + \epsilon)$ is a nonintegral eigenvalue of $\Delta$.
Then, by assumption, all of its conjugates lie in the interval $(m - \delta, m + \epsilon)$.
Since $\eta \in (m-2, m+2)$,
it follows by~\cite[Theorem~7.1]{BK99} that its conjugates form the set
$H := \left\{m + 2\cos\left({2\pi j \over k}\right) \;\middle|\; j \in \Z_k^*\right\}$
for some positive integer $k \not\in \{1, 2, 3, 4, 6\}$
(as for these values of $k$, the set $H$ may only contains integers).
Define
$$
\ell = \begin{cases}
{k-1 \over 2} & \text{if $k$ is odd,} \\
{k \over 2} - 1 & \text{if $k \equiv 0 \pmod{4}$, and} \\
{k \over 2} - 2 & \text{if $k \equiv 2 \pmod{4}$.}
\end{cases}
$$
Note that $\ell \in \Z_k^*$,
and ${1 \over 3} < {\ell \over k} < {1 \over 2}$
unless $k \in \{1, 2, 3, 4, 6, 10\}$,
implying $2\cos\left({2\pi \ell \over k}\right) < -1$ in this case.
Since $H \subset (m - \delta, m + \epsilon)$, it follows that $k = 10$,
i.e., $\eta = m + {1 \pm \sqrt{5} \over 2}$,
and then $\delta \ge {-1 + \sqrt{5} \over 2}$ and $\epsilon \ge {1 + \sqrt{5} \over 2}$
also follows.
Note also that $\eta = m+1 + 2\cos\left({2\pi \ell' \over 5}\right)$
for some $\ell' \in \Z_5^*$.
\end{proof}

When there are only finitely many candidates for eigenvalues,
a system of equations can be set to relate them to their multiplicities.

\begin{lemma}\label{lem:mult}
Let $\Delta$ be an $a_1$-regular graph with $k$ vertices
and distinct eigenvalues $\zeta_j$ $(1 \le j \le r)$
for some integer $r$.
Let $m_j$ be the multiplicity of the eigenvalue $\zeta_j$ $(1 \le j \le r)$.
Then, the following equations hold.

\begin{equation}\label{eqn:genmult}
\sum_{j=1}^r m_j = k \qquad
\sum_{j=1}^r m_j \zeta_j = 0 \qquad
\sum_{j=1}^r m_j \zeta_j^2 = k a_1
\end{equation}
\end{lemma}

\begin{proof}
Let $A$ be the adjacency matrix of $\Delta$.
The equations are derived by considering the
traces of the matrices $A^0, A^1, A^2$,
which are equal to the sums of their eigenvalues.
\end{proof}

Although the system of equations \eqref{eqn:genmult} only has three equations
and may therefore be underdetermined,
we are only interested in solutions
where the multiplicities are nonnegative integers.
We may also use additional information
-- for instance,
if the non-principal eigenvalues of the graph
are strictly smaller than the principal eigenvalue
(i.e., the graph is connected),
then the multiplicity of the principal eigenvalue is known to be $1$.
Furthermore, conjugate eigenvalues have the same multiplicity,
thus reducing the number of variables.

In the case the local graphs have at most four distinct eigenvalues,
they are walk-regular,
so we can also consider the numbers of closed walks through a vertex or an edge.

\begin{theorem}\label{thm:walks} {\rm \cite[\S3]{vD95}}
Let $\Delta$ be an $a_1$-regular graph on $k$ vertices
with at most four distinct eigenvalues
$a_1 = \zeta_1, \zeta_2, \zeta_3, \zeta_4$
of corresponding multiplicities $m_1, m_2, m_3, m_4$.
Then, for every positive integer $i$,
the number
$$
w_i := {1 \over k} \sum_{j=1}^4 m_j \zeta_j^i
$$
is the number of closed $i$-walks through a vertex of $\Delta$
and is therefore an integer;
if $i$ is odd, then $w_i$ is even.

Furthermore, the number
$$
w'_4 := {w_4 - a_1 (2a_1 - 1) \over 2}
$$
is the number of quadrangles through a vertex of $\Delta$
and is therefore an integer.
If the number of triangles through an edge of $\Delta$ is constant,
then there are precisely $2w'_4/a_1$ quadrangles through an edge of $\Delta$,
so this number is also an integer.
\qed
\end{theorem}

The above results constitute a new feasibility condition
for $Q$-polynomial distance-regular graphs.
We use the {\tt sage-drg} package~\cite{V18,Vdrg}
for the SageMath computer algebra system~\cite{Sage}
to compute these conditions.
In particular,
the package has been extended with methods
for computing Terwilliger polynomials,
further narrowing the possibilities for the eigenvalues of local graphs
following~\cite[Theorem~7.1]{BK99},
computing possible spectra of local graphs
using the integral solution finder (see~\cite[\S4.1]{GVW21}),
and verifying the integrality of the numbers of closed walks
in case the local graphs have at most four distinct eigenvalues.
This will allow us to derive some of the results presented in this paper.

We first show nonexistence
for some feasible intersection arrays of $Q$-polynomial distance-regular graphs.

\begin{theorem}\label{thm:nonex}
There is no distance-regular graph with intersection arrays
\begin{alignat*}{3}
\text{\rm (a)} &\quad& \{49, 40, 22; &\, 1, 5, 28\}
& \text{\rm ($750$ vertices)} &, \\
\text{\rm (b)} &\quad& \{109, 80, 22; &\, 1, 10, 88\}
& \text{\rm ($1\,200$ vertices)} &, \ \text{or} \\
\text{\rm (c)} &\quad& \{164, 121, 33; &\, 1, 11, 132\} \quad
& \text{\rm ($2\,420$ vertices)} &.
\end{alignat*}
\end{theorem}

\begin{proof}
We note that a distance-regular graph $\G$
with one of the intersection arrays above is $Q$-polynomial.
In all three cases,
the Terwilliger polynomial $T_2$ is nonnegative on the set of the form
$[\eta_{\min}, \gamma] \cup [\eta_{\max}, \overline{\gamma}]$,
where $\gamma$ and $\overline{\gamma}$ are conjugate algebraic integers.
Theorems~\ref{thm:bounds} and~\ref{thm:terwpoly} therefore imply
that the non-principal eigenvalues of $\G$ lie in the set
$[\eta_{\min}, \gamma] \cup \{\eta_{\max}\}$.

Let us give the values $\eta_{\min}, \gamma, \eta_{\max}, a_1$
for each intersection array.
\begin{alignat*}{6}
\text{\rm (a)} &\quad&
\eta_{\min} &= -{11 \over 3} &\qquad
\gamma &= {3 \over 2} - {\sqrt{345} \over 6} &&\approx -1.596 &\qquad
\eta_{\max} &= 3 &\qquad
a_1 &= 8 \\
\text{\rm (b)} &\quad&
\eta_{\min} &= -{11 \over 3} &\qquad
\gamma &= 9 - {8\sqrt{15} \over 3} &&\approx -1.328 &\qquad
\eta_{\max} &= 7 &\qquad
a_1 &= 28 \\
\text{\rm (c)} &\quad&
\eta_{\min} &= -{15 \over 4} &\qquad
\gamma &= {31 \over 2} - {11\sqrt{10} \over 2} &&\approx -1.893 &\qquad
\eta_{\max} &= 10&\qquad
a_1 &= 42
\end{alignat*}
In each case, we have $[\eta_{\min}, \gamma] \subset (-4, -1)$,
so the integral eigenvalues of the local graphs of $\G$
are all in the set $\{-3, -2, \eta_{\max}, a_1\}$.
In cases (a) and (c), we have $\gamma - (-3) < {1 + \sqrt{5} \over 2}$,
and Lemma~\ref{lem:noninteig} then implies
that there are no nonintegral eigenvalues;
for the case (b),
it implies that the only possible nonintegral eigenvalues
are ${-5 \pm \sqrt{5} \over 2}$,
both with the same multiplicity.
Also, in each case $\eta_{\max} < a_1$ holds,
so the principal eigenvalue $a_1$ has multiplicity $1$.

We can now use Lemma~\ref{lem:mult}
to establish the multiplicities of the eigenvalues.
However,
we find that in each case
the system of equations \eqref{eqn:genmult}
has no nonnegative integral solution.
We therefore conclude that a distance-regular graph
with one of the intersection arrays given above does not exist.
\end{proof}

\begin{remark}
The intersection array $\{49, 40, 22; 1, 5, 28\}$
is listed in~\cite{Bdrg,BCN} as a feasible parameter set.
Its formal dual has intersection array $\{112, 77, 16; 1, 16, 88\}$
and is ruled out by~\cite[Theorem~4.4.4]{BCN}.
The remaining two intersection arrays are formally self-dual.
All three intersection arrays are also listed as feasible parameters
for $Q$-polynomial association schemes in~\cite{WQpoly}.
\end{remark}

\section{Classical distance-regular graphs} \label{sec:clas}

A distance-regular graph $\G$ of diameter $d \ge 2$
is said to have classical parameters $\dbab$~\cite[\S6.1]{BCN}
if its intersection numbers satisfy
\begin{equation} \label{eq:clas}
b_i = ([d] - [i])(\beta - \alpha [i]), \qquad
c_i = [i](1 + \alpha [i-1]) \qquad
(0 \le i \le d),
\end{equation}
where $[i] = \sum_{j=0}^{i-1} b^j$.
Note that by~\cite[Proposition~6.2.1]{BCN},
the parameter $b$ is an integer distinct from $0$ and $-1$
whenever $d \ge 3$.
By~\cite[Proposition~4.3(i)]{GK15},
the Terwilliger polynomials $T_i$ ($2 \le i \le d-1$) of $\G$
have negative leading coefficients and roots
$$
\alpha b [d-1] - 1, \quad \beta - \alpha - 1, \quad -1, \quad  -b-1.
$$
Since distance-regular graphs with classical parameters
are $Q$-polynomial~\cite[Corollary~8.4.2]{BCN},
we may apply results from Section~\ref{sec:spectra}.

\begin{lemma}\label{lem:small intervals}
Suppose that $b \ge 1$ and $\beta \le \alpha [d]$ hold.
Then the eigenvalues of the local graphs of $\G$ are restricted to
\begin{equation*}\label{eq:small intervals}
[-b-1, -1] \cup \{\beta - \alpha - 1\} \cup \{a_1\}.
\end{equation*}
\end{lemma}
\begin{proof}
Since $b \ge 1$ holds,
$\beta > \alpha$ follows from $b_1 > 0$.
We then have $-b-1 < -1 < \beta - \alpha - 1 \le \alpha b [d-1] - 1$,
and the polynomials $T_i$ are nonnegative on
$$
[-b-1, -1] \cup [\beta - \alpha - 1, \alpha b [d-1] - 1] .
$$
Since, in the notation of Theorem~\ref{thm:bounds}, $\eta_{\min} = -b - 1$, $\eta_{\max} = \beta - \alpha - 1$
(see~\cite[Corollary~8.4.2]{BCN})
and $a_1 = k - b_1 - 1 = \beta + \alpha b [d-1] - 1$,
the result follows by Theorems~\ref{thm:bounds} and~\ref{thm:terwpoly}.
\end{proof}

We first derive nonexistence for two feasible intersection arrays
of classical distance-regular graphs with $d = b = 3$.

\begin{theorem}\label{thm:nonexb3}
There is no distance-regular graph with classical parameters
$(3, 3, 1/2, 5)$ or $(3, 3, 5/4, 8)$,
i.e., with intersection arrays
\begin{alignat*}{2}
\{65, 54, 27; &\, 1, 6, 39\} &\quad \text{\rm ($1\,056$ vertices)} & \ \text{or} \\
\{104, 81, 27; &\, 1, 9, 78\} &\quad \text{\rm ($1\,365$ vertices)}
&.
\end{alignat*}
\end{theorem}

\begin{proof}
Let $\G$ be a distance-regular graph
with one of the intersection arrays given above.
In both cases,
the value $\eta_{\max} = \beta - \alpha - 1$ is a nonintegral rational number,
and by Lemma~\ref{lem:small intervals},
the non-principal eigenvalues of the local graphs are thus all negative.
Therefore, the local graphs of $\G$ are unions of cliques of size $a_1 + 1 < k$.
However,
Theorem~\ref{thm:bounds} implies
that the local graphs of $\G$ should be connected,
which leads to a contradiction.
We therefore conclude that a distance-regular graph
with one of the intersection arrays given above does not exist.
\end{proof}

\begin{remark}
Both of the above intersection arrays are listed as feasible parameters
for $Q$-polynomial association schemes in~\cite{WQpoly}.
\end{remark}

\underline{From now on we focus on the case $b = 2$.}
Under the hypothesis of Lemma~\ref{eq:small intervals},
it follows from  Lemma~\ref{lem:noninteig} that
the spectrum of the local graphs of $\G$ is restricted to
$$
\left\{
{[\beta + \alpha(2^d - 2) - 1]}^1, \quad
{[\beta - \alpha - 1]}^{m_1}, \quad
{[-1]}^{m_2}, \quad
{[-2]}^{m_3}, \quad
{[-3]}^{m_4}
\right\}
$$
for some nonnegative integers $m_1, m_2, m_3, m_4$.
Following Lemma~\ref{lem:mult},
these are further constrained by the following system of equations.
\begin{equation}\label{eqn:mult}
\begin{alignedat}{4}
m_1 &\,+\,& m_2 &\,+\,& m_3 &\,+\,& m_4 &= (2^d - 1) \beta - 1 \\
(\beta - \alpha - 1) m_1 &\,-\,& m_2 &\,-\,& 2 m_3 &\,-\,& 3 m_4
&= -\beta - \alpha(2^d - 2) + 1 \\
(\beta - \alpha - 1)^2 m_1 &\,+\,& m_2 &\,+\,& 4 m_3 &\,+\,& 9 m_4
&= ((2^d - 2) (\beta-\alpha) + 1) (\beta + \alpha(2^d - 2) - 1)
\end{alignedat}
\end{equation}

We can now derive nonexistence for some more feasible parameter sets.

\begin{theorem}\label{thm:nonex2}
There is no distance-regular graph with classical parameters
$(3, 2, 2, 9)$, $(3, 2, 5, 21)$ or $(6, 2, 2,\allowbreak 107)$,
i.e., with intersection arrays
\begin{alignat*}{2}
\{63, 42, 12; &\, 1, 9, 49\} & \text{\rm ($430$ vertices)} &, \\
\{147, 96, 24; &\, 1, 18, 112\} & \text{\rm ($1\,100$ vertices)}
&, \ \text{or} \\
\{6741, 6510, 6060, 5208, 3696, 1440; &\, 1, 9, 49, 225, 961, 3969\} \quad
& \text{\rm ($87\,725\,820\,468$ vertices)} &.
\end{alignat*}
\end{theorem}

\begin{proof}
In each of the given cases,
the system of equations \eqref{eqn:mult}
has no nonnegative integral solution.
\end{proof}

\begin{remark}
The intersection array $\{63, 42, 12; 1, 9, 49\}$
is listed in~\cite{Bdrg,BCN} as a feasible parameter set,
while the intersection array $\{147, 96, 24; 1, 18, 112\}$
is mentioned in~\cite{JV17} as a feasible example
of classical parameters with $\alpha \ne b-1, b$.
Both are also listed as feasible parameters
for $Q$-polynomial association schemes in~\cite{WQpoly}.
Note that Koolen et al.~\cite{KOOLEN2026104378} recently showed that
$\alpha\leq 2$ if $b=2$ and $D\geq 12$.
\end{remark}

\underline{For the rest of the paper we assume $\alpha = 1$.}
The system of equations \eqref{eqn:mult} now simplifies as follows.
\begin{equation}\label{eqn:mult1}
\begin{alignedat}{4}
m_1 &\,+\,& m_2 &\,+\,& m_3 &\,+\,& m_4 &= (2^d - 1) \beta - 1 \\
(\beta - 2) m_1 &\,-\,& m_2 &\,-\,& 2 m_3 &\,-\,& 3 m_4
&= -2^d - \beta + 3 \\
(\beta - 2)^2 m_1 &\,+\,& m_2 &\,+\,& 4 m_3 &\,+\,& 9 m_4
&= ((2^d - 2) (\beta-1) + 1) (2^d + \beta - 3)
\end{alignedat}
\end{equation}
We can express its general solution in terms of a parameter $t$.
\begin{equation}\label{eqn:sol1}
\begin{aligned}
m_1 &= {(2^d - 2)(2^d + \beta - 1) - t \over \beta} \\
m_2 &= t \cdot {\beta + 1 \over 2} - (2^d - \beta - 1)(2^d - 1) \\
m_3 &= ((2^d - 2)(2^d - 1) - t (\beta + 1)) \cdot {\beta - 1 \over \beta} \\
m_4 &= t \cdot {\beta - 1 \over 2}
\end{aligned}
\end{equation}
As $\beta \ge 2^{d-1}$ holds by~\cite[Proposition~1(iii)]{JV17},
we must have
\begin{equation}\label{eqn:interval}
{2(2^d - \beta - 1)(2^d - 1) \over \beta + 1}
\le t \le {(2^d - 2)(2^d - 1) \over \beta + 1}
\end{equation}
for the solution to be nonnegative.
As $\beta$ is integral by~\cite[Lemma~3]{JV17},
$t$ must be an integer such that
\begin{equation}\label{eqn:modulo}
t \equiv (2^d - 2)(2^d - 1) \pmod{\beta}
\end{equation}
for the solution to be integral.
In particular, if $\beta$ is even, $t$ must also be even.

Let us now find the parameter sets with $d \le 8$
for which there is no nonnegative integral solution
to the system of equations \eqref{eqn:mult1}.

\begin{proposition}\label{prop:nonex1}
There is no distance-regular graph
with classical parameters $(d, 2, 1, \beta)$
such that
\begin{itemize}
\item $d = 4$, $\beta \in \{8, 10, 12\}$,
\item $d = 5$, $\beta \in \{16, 17, 19, 20, 21, 28\}$,
\item $d = 6$, $\beta \in \{32, 33, 34, 35, 36, 38, 40, 46, 49, 54, 60\}$,
\item $d = 7$,
$\beta \in \{64, 65, 66, 67, 69, 70, 71, 72, 73, 74, 77, 79, 81, 84, 85, 92,
99, 124\}$,
\item $d = 8$,
$\beta \in \{128, 129, 130, 131, 133, 134, 135, 136, 137, 139, 140, 141, 151,
152, 155, 158,\allowbreak 160, 165, 168,\allowbreak
174, 175, 183, 184, 190, 202, 238, 252\}$.
\end{itemize}
\end{proposition}

\begin{proof}
In each of the given cases,
there is no integer $t$ satisfying
\eqref{eqn:interval} and \eqref{eqn:modulo}.
\end{proof}

\begin{remark}
The cases with $d = 4$ and $\beta = 8, 10$
and intersection arrays
\begin{alignat*}{2}
\{120, 98, 60, 8; &\, 1, 6, 28, 120\} &\quad \text{\rm ($6\,561$ vertices)}
&\ \text{and} \\
\{150, 126, 84, 24; &\, 1, 6, 28, 120\} &\quad \text{\rm ($14\,641$ vertices)}
\end{alignat*}
are listed in~\cite{Bdrg} as feasible parameter sets.
The case with $(d, \beta) = (5, 16)$ is not feasible
since it has intersection number $p^5_{55} < 0$.
\end{remark}

We note that the cases $\beta = 2^{d-1}$ and $\beta = 2^d - 4$
appear as nonexistent for all the values of $d$
given in Proposition~\ref{prop:nonex1}.
Nonexistence of the family of distance-regular graphs
with classical parameters $(d, 2, 1, 2^{d-1})$ and $d \ge 3$
follows from~\cite[Theorem~8]{JV17}.
We now also show nonexistence of the infinite family for the latter case.

\begin{theorem}\label{thm:nonex4}
There is no distance-regular graph
with classical parameters $(d, 2, 1, 2^d - 4)$ and $d \ge 3$.
\end{theorem}

\begin{proof}
We claim that a nonnegative solution \eqref{eqn:sol1}
of the system of equations \eqref{eqn:mult1}
with $\beta = 2^d - 4$ satisfies
$$
6 < 6 \cdot {2^d - 1 \over 2^d - 3} \le t
\le {(2^d - 2)(2^d - 1) \over 2^d - 3} < 2^d + 2 .
$$
The inner two inequalities follow by substituting into \eqref{eqn:interval},
while the leftmost inequality is obvious.
We can verify that the rightmost inequality simplifies to $2^{d+1} > 8$,
which is true for all $d \ge 3$.

By \eqref{eqn:modulo}, we have $t \equiv 6 \pmod{2^d - 4}$.
However, we see that no such value of $t$ satisfies the above condition.
Therefore,
a distance-regular graph with classical parameters $(d, 2, 1, 2^d - 4)$
and $d \ge 3$ does not exist.
\end{proof}

We now consider the cases when the system of equations \eqref{eqn:mult1}
has a nonnegative integral solution.
Let us first show that such a solution must be unique.

\begin{lemma}\label{lem:unique}
Given integers $d \ge 3$ and $\beta \le 2^d - 1$,
the system of equations \eqref{eqn:mult1}
has at most one nonnegative integral solution.
\end{lemma}

\begin{proof}
Let us show that the nonnegative solutions \eqref{eqn:sol1}
of the system of equations \eqref{eqn:mult1}
have $t$ from an interval of length less than $\beta$.
Using \eqref{eqn:interval}, let us show the following inequality:
\begin{align*}
{(2^d - 2)(2^d - 1) \over \beta + 1} -
{2(2^d - \beta - 1)(2^d - 1) \over \beta + 1} &< \beta\text{, i.e.,} \\
(2\beta - 2^d)(2^d - 1) &< \beta (\beta + 1).
\intertext{We now write $\beta = 2^d - r$:}
(2^d - 2r)(2^d - 1) &< (2^d - r)(2^d - r + 1)\\
r(3 - r) &< 2^{d+1}.
\end{align*}
The above is true for all $r$ and $d \ge 3$.
By \eqref{eqn:modulo} it then follows that a solution \eqref{eqn:sol1}
is nonnegative and integral for at most one choice of $t$.
\end{proof}

\begin{corollary}\label{cor:unique}
Let $\G$ be a distance-regular graph with classical parameters
$(d, 2, 1, \beta)$ such that $d \ge 3$ and $\beta \le 2^d - 1$.
Then all its local graphs have the same spectrum. \qed
\end{corollary}

Let us now consider when local graphs of $\G$
have at most four distinct eigenvalues
(i.e., $m_i = 0$ for some $i$).
If $t$ attains the lower bound of \eqref{eqn:interval},
then $m_2 = 0$.
For this to occur, $\beta + 1$ must divide $2^{d+1}(2^d - 1)$,
and the following congruence must also hold.
\begin{alignat*}{2}
t &\equiv (2^d - 2)(2^d - 1) &&\pmod{\beta} \\
{2(2^d - \beta - 1)(2^d - 1) \over \beta + 1}
&\equiv (2^d - 2)(2^d - 1) \cdot {\beta + 1 \over \beta + 1} &&\pmod{\beta} \\
{2^d (2^d - 1) \over \beta + 1} &\equiv 0 &&\pmod{\beta} \\
2^d (2^d - 1) &\equiv 0 &&\pmod{\beta (\beta + 1)}
\end{alignat*}
Since $\beta \ge 2^{d-1}$,
the above congruence can only hold when $\beta = 2^d - 1$,
in which case we have $t = 0$ and $m_4 = 0$,
leading to the characterization of the bilinear forms graphs in~\cite{GK19}.

If, on the other hand, $t$ attains the upper bound of \eqref{eqn:interval},
then $m_3 = 0$ and the local graphs of $\G$ have four distinct eigenvalues.
For this to occur, $\beta + 1$ must divide $(2^d - 2)(2^d - 1)$,
and the following congruence must also hold.
\begin{alignat*}{2}
t &\equiv (2^d - 2)(2^d - 1) &&\pmod{\beta} \\
{(2^d - 2)(2^d - 1) \over \beta + 1}
&\equiv (2^d - 2)(2^d - 1) \cdot {\beta + 1 \over \beta + 1} &&\pmod{\beta} \\
{(2^d - 2)(2^d - 1) \over \beta + 1} \cdot \beta &\equiv 0 &&\pmod{\beta}
\end{alignat*}
We see that the above congruence does not impose a further restriction.

From now on we will thus assume that $\beta + 1$ divides $(2^d - 2)(2^d - 1)$.
Considering Corollary~\ref{cor:unique},
it follows that all local graphs of $\G$ have spectrum
$$
\left\{
{[\beta +  2^d - 3]}^1, \quad
{[\beta - 2]}^{m_1}, \quad
{[-1]}^{m_2}, \quad
{[-3]}^{m_4}
\right\},
$$
where the multiplicities of the non-principal eigenvalues are as follows.
\begin{align*}
m_1 &= {(2^d - 2)(2^d + \beta) \over \beta + 1} \\
m_2 &= (\beta - 2^{d-1})(2^d - 1) \\
m_4 &= (2^{d-1} - 1)(2^d - 1) \cdot {\beta - 1 \over \beta + 1}
\end{align*}

We may now further constrain $\beta$.

\begin{theorem}\label{thm:nonext}
Let $\G$ be a distance-regular graph with classical parameters
$(d, 2, 1, \beta)$ such that $d \ge 3$,
$\beta \le 2^d - 1$ and $(\beta + 1) \,|\, (2^d - 2)(2^d - 1)$.
Then $\beta = 2^d - 2$.
\end{theorem}

\begin{proof}
By Theorem~\ref{thm:walks},
the value
$$
w_i := {1 \over \beta (2^d - 1)} \left(
(\beta + 2^d - 3)^i + m_1 (\beta - 2)^i + m_2 (-1)^i + m_4 (-3)^i
\right)
$$
counting closed walks of length $i$ through each vertex
should be integral for all $i \ge 3$, and even whenever $i$ is odd.
Let us express $w_3$.
$$
w_3 = \beta^2 + (2^d - 5) \beta + 10 - 2^{d+2} + {2^d (2^d - 2) \over \beta}
$$
Therefore, $\beta$ should divide $2^d (2^d - 2)$
-- in this case, $w_3$ is even since $\beta < 2^d$ holds.
It then follows that there are positive integers $r, s$ such that
\begin{align}
\beta r &= 2^d (2^d - 2) = 2^{d+1} (2^{d-1} - 1) \quad \text{and} \notag \\
(\beta + 1) s &= (2^d - 1)(2^d - 2) . \label{eqn:s}
\end{align}
If $\beta$ were odd, then it would divide $2^{d-1} - 1$,
which contradicts $\beta \ge 2^{d-1}$.
Therefore, $\beta$ is even.

Let us set $u := r - s$.
We now derive
\begin{align}
(\beta + 1) u &= r + 2^d - 2 \qquad \text{and} \label{eqn:u} \\
u &= {(2^d - 2)(2^d + \beta) \over \beta (\beta + 1)}. \notag
\end{align}
Since $\beta \ge 2^{d-1}$, we obtain
$$
u \le 6 \cdot {2^d - 2 \over 2^d + 2}
$$
and therefore $u \le 5$.
On the other hand, since $\beta \le 2^d - 1$,
we have $r \equiv 0 \pmod{4}$
and thus $r + 2^d - 2 \equiv 2 \pmod{4}$.
As $\beta + 1$ is odd, we then also have $u \equiv 2 \pmod{4}$,
from which $u = 2$ and $r = s + 2$ follow.

By eliminating $\beta$ from equations (\ref{eqn:s}--\ref{eqn:u}),
we obtain
$$
2 (2^d - 1)(2^d - 2) = s (2^d + s) .
$$
For all $d \ge 3$,
the only positive solution of the above equation is $s = 2^d - 2$ ,
which then gives us $\beta = 2^d - 2$.
\end{proof}

\begin{remark} \label{rem:walks}
For distance-regular graphs with classical parameters $(d, 2, 1, 2^d - 2)$,
we get no new restrictions
from considering the number of walks of lengths greater than $3$
-- i.e., it can be shown that the numbers $w_i$ $(i > 3)$ are integers,
and that they are even whenever $i$ is odd.
\end{remark}

We now list the parameter sets with $d \le 8$
which are ruled out by Theorem~\ref{thm:nonext}.

\begin{corollary}\label{cor:nonext}
There is no distance-regular graph
with classical parameters $(d, 2, 1, \beta)$ such that
$(d, \beta) \in \{(3, 5), (4, 9), (4, 13), (5, 29), (6, 41), (6, 61),
(7, 125), (8, 169), (8, 253)\}$. \qed
\end{corollary}

\begin{remark}
The cases with $(d, \beta) = (3, 5), (4, 9)$
and intersection arrays
\begin{alignat*}{2}
\{35, 24, 8; &\, 1, 6, 28\} &\quad \text{\rm ($216$ vertices)}
&\ \text{and} \\
\{135, 112, 84, 24; &\, 1, 6, 28, 120\} &\quad \text{\rm ($10\,000$ vertices)}
\end{alignat*}
are listed in~\cite{Bdrg} as feasible parameter sets.
The former is also listed in~\cite{BCN}
and was shown to be infeasible in~\cite[Theorem~3]{JV12}
as a part of an infinite family.
An alternative proof of nonexistence,
of which the above proof is a generalization,
was given in~\cite[Theorem~5.1]{GK19}.
\end{remark}

We note that the cases given in Corollary~\ref{cor:nonext}
can be covered by two infinite families of parameter sets.
Let us show nonexistence for both families.

\begin{theorem}\label{thm:nonex3}
There is no distance-regular graph
with classical parameters $(d, 2, 1, 2^d - 3)$ and diameter $d \ge 3$
or $\left(d, 2, 1, {2^{d+1} - 5 \over 3}\right)$ and even diameter $d \ge 4$.
\end{theorem}

\begin{proof}
In either case, $\beta + 1$ divides $(2^d - 2)(2^d - 1)$.
Since we have $\beta \ne 2^d - 2$,
Theorem~\ref{thm:nonext} implies nonexistence.
\end{proof}

Let us now summarize the results regarding the spectra of local graphs
of distance-regular graphs with classical parameters $(d, 2, 1, \beta)$.

\begin{corollary} \label{cor:loceig}
Let $\G$ be a distance-regular graph
with classical parameters $(d, 2, 1, \beta)$,
where $d \ge 3$ and $\beta \le 2^d - 1$.
Then one of the following holds.
\begin{itemize}
\item $2^{d-1} + 1 \le \beta \le 2^d - 5$,
and the local graphs of $\G$ are all cospectral with $5$ distinct eigenvalues.
\item $\beta = 2^d - 2$,
and the local graphs of $\G$ are all cospectral with $4$ distinct eigenvalues.
\item $\beta = 2^d - 1$ and $\G \cong H_2(d, d)$
(in particular,
the local graphs of $\G$ are all isomorphic to $H(2, 2^d - 1)$,
which has $3$ distinct eigenvalues).
\end{itemize}
\end{corollary}

\begin{proof}
By Corollary~\ref{cor:unique}, all local graphs of $\G$ are cospectral.
The arguments above show that
the local graphs of $\G$ have three distinct eigenvalues if $\beta = 2^d - 1$
(in which case, $\G \cong H_2(d, d)$ by~\cite[Theorem~1.3]{GK19}),
four distinct eigenvalues if $\beta + 1$ divides $(2^d - 2)(2^d - 1)$,
which implies $\beta = 2^d - 2$ by Theorem~\ref{thm:nonext},
and five distinct eigenvalues otherwise.
Assuming the latter,
note that~\cite[Proposition~1(iii), Theorem~8]{JV17}
imply that $\beta \ge 2^{d-1} + 1$,
while Theorems~\ref{thm:nonex4} and~\ref{thm:nonex3} imply $\beta \le 2^d - 5$.
\end{proof}

%%%%%%%%%%%%%%%%%%%%%%%%%%%%%%%%%%%%%%%%%%%%%%%%%%%%%%%%%%%%%%%%%%%%%%%%
\section{Local graphs in classical DRGs
with $(d, b, \alpha, \beta) = (d, 2, 1, 2^d-2)$}

We will now investigate what local graphs
a distance-regular graphs $\G$ with classical parameters $(d, 2, 1, 2^d - 2)$
and $d \ge 3$ might have.
By the results in Section~\ref{sec:clas},
the local graphs of $\G$ have spectrum equal to
\begin{equation}\label{eq:classpectrum}
\left\{
{[2^{d+1} - 5]}^1, \quad
{[2^d - 4]}^{2^{d+1} - 4}, \quad
{[-1]}^{(2^{d-1} - 2)(2^d - 1)}, \quad
{[-3]}^{(2^{d-1} - 1)(2^d - 3)}
\right\}.
\end{equation}

Let $r,t$ be positive integers.
Let $\Theta_t$ denote a graph on $2t+1$ vertices which is the disjoint union of two $t$-cliques and an isolated vertex.
Let $\Lambda_r$ denote a graph on the vertex set $\{(i,j)\mid 1\leq i,j\leq r, i\ne j\}$
with two distinct vertices $(i,j)$, $(i',j')$ being adjacent if
$i=i'$ or $j=j'$ or $(i, j) = (j', i')$.

\begin{lemma}\label{lm:localTheta}
Let $t \ge 1$.
A locally $\Theta_{t}$ graph with $(t+1)(t+2)$ vertices
is isomorphic to $\Lambda_{t+2}$.
\end{lemma}

\begin{proof}
Let $\Delta$ be a locally $\Theta_t$ graph with $(t+1)(t+2)$ vertices.
Since $\Theta_t$ is a union of cliques,
it follows that each pair of distinct maximal cliques of $\Delta$
has at most one common vertex
-- in other words,
a vertex $u$ outside of a maximal clique $C$
has at most one neighbour in $C$.
Furthermore, the maximal cliques of $\Delta$ have sizes $2$ and $t+1$.
If $t = 1$,
then $\Delta$ is a triangle-free cubic graph on $6$ vertices
and is therefore isomorphic to $K_{3,3} \cong \Lambda_3$.
For the rest of the proof, we thus assume $t \ge 2$.

Let $A_1 := \{(1, j) \mid 2 \le j \le t+2\}$ be a $(t+1)$-clique of $\Delta$,
and let $B_j$ ($2 \le j \le t+2$)
be the $(t+1)$-clique distinct from $A_1$ containing $(1, j)$.
By the argument above, these cliques are pairwise disjoint.
Furthermore,
let $B_1 := \{(i, 1) \mid 2 \le i \le t+2\}$,
where $(i, 1)$ is the neighbour of $(1, i)$ ($2 \le i \le t+1$)
not contained in any of the $(t+1)$-cliques $B_j$ ($2 \le j \le t+2$),
and note that $|B_1| = t+1$ (i.e., $(i, 1) = (i', 1)$ iff $i = i'$).
As $\left|\bigcup_{j=1}^{t+2} B_j\right| = (t+1)(t+2)$,
each vertex of $\Delta$ is contained
in precisely one of the sets $B_j$ ($1 \le j \le t+2$).

Since each of the vertices $(i, 1)$
has at most one neighbour in $B_j$ ($2 \le i, j \le t+2$),
it follows that $B_1$ is also a $(t+1)$-clique.
Let $A_i$ ($2 \le i \le t+2$)
be the $(t+1)$-clique distinct from $B_1$ containing $(i, 1)$.
Similarly as before,
we see that each vertex of $\Delta$ is contained
in precisely one of the $(t+1)$-cliques $A_i$ ($1 \le i \le t+2$).

Since the vertex $(i, 1)$ ($2 \le i \le t+2$)
is adjacent to the vertex $(1, i) \in B_i \setminus A_i$,
the set $A_i \cap B_j$ ($2 \le j \le t+2$) is empty when $i = j$,
and contains precisely one vertex otherwise
-- call this vertex $(i, j)$,
and let $(i', j')$ be the neighbour of $(i, j)$ outside $A_i \cup B_j$.
Clearly, we must have $(i', j') \in A_j \cap B_i$,
and therefore $(i, j)$ is adjacent to $(j, i)$.
This establishes that the graph $\Delta$ is isomorphic to $\Lambda_{t+2}$.
\end{proof}

\begin{lemma}\label{lm:LambdaSpectrum}
Let $t \ge 1$.
The graph $\Lambda_{t+2}$ has $k_t:=t^2+3t+2$ vertices and spectrum $\{{[\zeta_i]}^{m_i}\}_{i=0}^3$ equal to
\begin{equation}\label{eq:spectrum}
\left\{
{[2t+1]}^1, \quad
{[t-1]}^{2(t+1)}, \quad
{[-1]}^{(t-1)(t+2)/2}, \quad
{[-3]}^{t(t+1)/2}
\right\}.
\end{equation}
\end{lemma}

\begin{proof}
Let $X = \{(i, j) \mid 1 \le i, j \le t+2, i \ne j\}$
be the vertex set of the graph $\Lambda_{t+2}$.
Clearly, we have $|X| = (t+1)(t+2) = k_t$.
Denote the multiplicities of the eigenvalues $2t+1$, $t-1$, $-1$ and $-3$
of the graph $\Lambda_{t+2}$ by $m_0, m_1, m_2, m_3$,
respectively.
As $\Lambda_{t+2}$ is a $(2t+1)$-regular connected graph,
we have $m_0 = 1$.
To determine the multiplicities of the remaining eigenvalues,
we will construct an appropriate number of linearly independent eigenvectors.

Let $C$ be a $(t+1)$-clique of $\Lambda_{t+2}$.
Define the vector $u^{(C)}$ by
$$
u^{(C)}_{x} = \begin{cases}
t+1 & \text{if $x \in C$, and} \\
-1 & \text{otherwise.}
\end{cases}
\qquad (x \in X)
$$
Let $A_i = \{(i, j) \mid 1 \le j \le t+2, i \ne j\}$
and $B_i = \{(j, i) \mid 1 \le j \le t+2, i \ne j\}$ ($1 \le i \le t+2$)
be $(t+1)$-cliques of $\Lambda_{t+2}$.
Then $u^{(A_i)}$ and $u^{(B_i)}$ ($1 \le i \le t+2$)
are eigenvectors of $\Lambda_{t+2}$ for the eigenvalue $t-1$.
It can be easily checked that the set
$\{u^{(A_i)}, u^{(B_i)} \mid 1 \le i \le t+1\}$
is linearly independent.
Therefore, we have $m_1 \ge 2(t+1)$.

Next,
we define the vectors $v^{\{(i, j), (k, \ell)\}}$
for distinct $1 \le i, j, k, \ell \le t+2$ by
$$
v^{\{(i, j), (k, \ell)\}}_x = \begin{cases}
1 & \text{if $x \in \{(i, j), (j, i), (k, \ell), (\ell, k)\}$,} \\
-1 & \text{if $x \in \{(i, \ell), (j, k), (k, j), (\ell, i)\}$, and} \\
0 & \text{otherwise.}
\end{cases}
\qquad (x \in X)
$$
These vectors are eigenvectors of $\Lambda_{t+2}$ for the eigenvalue $-1$.
Define the sets
\begin{align*}
V_1 &= \{v^{\{(1, 2), (i, j)\}} \mid 3 \le i < j \le t+2\} \quad \text{and} \\
V_2 &= \{v^{\{(1, 3), (2, i)\}} \mid 4 \le i \le t+2\}.
\end{align*}
Since $v^{\{(1, 2), (i, j)\}}$ is the unique vector of $V_1$
such that its support contains the vertex $(i, j)$ ($3 \le i < j \le t+2$)
and there is no such vector in $V_2$,
and $v^{\{(1, 3), (2, i)\}}$ is the unique vector of $V_2$
such that its support contains the vertex $(2, i)$ ($4 \le i \le t+2$),
it follows that $V_1 \cup V_2$ is linearly independent.
Therefore, we have $m_2 \ge {t(t-1) \over 2} + t-1 = {(t-1)(t+2) \over 2}$.

Finally,
we define vectors $w^{(i, j, k)}$ for distinct $1 \le i, j, k \le t+2$ by
$$
w^{(i, j, k)}_x = \begin{cases}
1 & \text{if $x \in \{(i, j), (j, k), (k, i)\}$,} \\
-1 & \text{if $x \in \{(i, k), (j, i), (k, j)\}$, and} \\
0 & \text{otherwise.}
\end{cases}
\qquad (x \in X)
$$
These vectors are eigenvectors of $\Lambda_{t+2}$ for the eigenvalue $-3$.
The set $W = \{w^{(1, i, j)} \mid 2 \le i < j \le t+2\}$
is linearly independent,
as $w^{(1, i, j)}$ is the unique vector of $W$
such that its support contains the vertex $(i, j)$ ($2 \le i < j \le t+2$).
We thus have $m_3 \ge {t(t+1) \over 2}$.

We conclude by noting that $m_0 + m_1 + m_2 + m_3 \ge t^2 + 3t + 2 = k_t$.
All the lower bounds for the multiplicities are therefore met with equality,
so the graph $\Lambda_{t+2}$ has the specified spectrum.
\end{proof}

The following is the main result of this section. Our proof relies on the technique developed in~\cite{grids,MR4825202,thinGrassmann}.

\begin{proposition}\label{pr:delta}
Let $\Delta$ be a graph with spectrum \eqref{eq:spectrum}, where $t\geq 292$, such that each pair of nonadjacent vertices of $\Delta$
has at most $5$ common neighbours.
Then $\Delta$ is
isomorphic to $\Lambda_{t+2}$.
\end{proposition}

In what follows,
let $\Delta$ be a graph satisfying the hypothesis of Proposition~\ref{pr:delta}.
We call a maximal clique of $\Delta$ a \textbf{line}
if it contains at least $\kappa t+1$ vertices,
where $\kappa:=0.67124$.
We first show in Lemma~\ref{lm:2lines}
that every vertex of $\Delta$ lies in exactly two lines.
In Lemma~\ref{lm:linesize},
it is proved that every line has the same size $t+1$
and every two distinct intersecting lines share precisely one vertex.

\begin{lemma}\label{lm:pre}
\begin{itemize}
\item[$(1)$] A clique of $\Delta$ has size at most $t+1$;
if equality is attained,
then every vertex outside of the clique
is adjacent to exactly one of its vertices.
\item[$(2)$] For each vertex of $\Delta$,
its local graph contains a coclique of size at most $9$.
\item[$(3)$] The number of edges in a local graph of $\Delta$ equals $t(t-1)$.
\end{itemize}
\end{lemma}
\begin{proof}
Let $L$ be a clique of size $\ell$ in $\Delta$.
The partition $\{L,\Delta\setminus L\}$ of the vertex set of $\Delta$
has quotient matrix
\[
  Q:=\kbordermatrix{
    & L & \Delta\setminus L \\
  L & \ell-1 & 2t+2-\ell \\
 \Delta\setminus L & \frac{(2t+2-\ell)\ell}{k_t-\ell} & 2t+1-\frac{(2t+2-\ell)\ell}{k_t-\ell}
  }
\]
with eigenvalues $2t+1$ and $\frac{k_t - \ell (1 + t + t^2)}{\ell - k_t}$.
By the eigenvalue interlacing theorem (see~\cite[Section~2.5]{BrouwerHaemers2012}),
the second largest eigenvalue of $Q$
is at most the second largest eigenvalue of $\Delta$,
which is $t-1$.
This gives $\ell\leq t+1$.
In case of equality,
the partition $\{L,\Delta\setminus L\}$ is equitable;
thus, $Q_{\Delta\setminus L, L}=\frac{(2t+2-\ell)\ell}{t^2+3t+2-\ell}=1$
is the number of neighbors in $L$ of a vertex from $\Delta\setminus L$.
This proves (1).

The smallest eigenvalue of a complete bipartite graph
with parts of size $1$ and $c$ is $-\sqrt{c}$.
Thus, by the eigenvalue interlacing theorem,
a local graph of $\Delta$ may have a $c$-coclique
only if $-\sqrt{c}\geq -3$.
This shows (2).

Since $\Delta$ has just $4$ distinct eigenvalues,
it follows from Theorem~\ref{thm:walks}
that the number of triangles through a vertex of $\Delta$ equals
\begin{align*}
{w_3 \over 2} = \frac{1}{2k_t}\sum_{i=0}^{3}m_i\zeta_i^3 = t(t-1),
\end{align*}
which implies (3),
and the lemma follows.
\end{proof}

\begin{lemma}\label{lm:2lines}
Every vertex of $\Delta$ belongs to exactly two lines.
\end{lemma}
\begin{proof}
Fix a vertex $\infty$ of $\Delta$.
Let $C$ be a maximum coclique of $\Delta(\infty)$ with vertex set $\{x_1,x_2,\ldots,x_c\}$.
By Lemma~\ref{lm:pre}(2), one has $c\leq 9$.
Define
\begin{align*}
P &:= \{y\in \Delta(\infty) \mid
\text{$y$ has at least two neighbours in $C$}\}, \\
U_i &:= \{x_i\}\cup \{y\in \Delta(\infty) \mid
\text{$y$ has only $x_i$ as its neighbour in $C$}\} \quad (1\leq i\leq c).
\end{align*}

The maximality of $C$ implies that
$\{P,U_1,\ldots,U_c\}$ is a partition of the vertex set of $\Delta(\infty)$
and each $U_i$ induces a complete subgraph in $\Delta(\infty)$.
Put $p:=|P|$, and $u_i:=|U_i|$ for $1\leq i\leq c$.
Since, by assumption, a pair of nonadjacent vertices of $\Delta$ has at most $5$ common neighbours,
the vertices $x_i$ and $x_j$ with $i\ne j$
have at most $4$ common neighbours in $\Delta(\infty)$;
thus, we obtain
\begin{equation}\label{eq:p}
p\leq 4{c\choose 2}
\end{equation}

Let us estimate the number of edges in $\Delta(\infty)$:
\begin{itemize}
\item \underline{Edges in $U_i$}.
Since $U_i$ is a clique, it contains $u_i(u_i-1)/2$ edges.
\item \underline{Edges in $P$}.
$P$ contains at most $p(p-1)/2$ edges.
\item \underline{Edges between $C$ and $P$}.
There are at most $cp$ edges between $C$ and $P$.
\item \underline{Edges between $U_i\setminus\{x_i\}$
and $\bigcup_{j\ne i}U_j\cup P$}.
Each vertex of $U_i\setminus\{x_i\}$
has at most $4$ common neighbors with $x_j$, $j\ne i$;
thus, there are at most $4(u_i-1)(c-1)$ edges
from $U_i\setminus\{x_i\}$ to $\bigcup_{j\ne i}U_j\cup P$.
\end{itemize}

On the other hand, by Lemma~\ref{lm:pre}(3),
the number of edges in $\Delta(\infty)$ equals $t(t-1)$.
Hence, we obtain
\[
2t(t-1)\leq \sum_{i=1}^c u_i(u_i-1) + p(p-1) + 2cp + \sum_{i=1}^c 4(c-1)(u_i-1),
\]
which, by using Equation \eqref{eq:p} and $\sum_{i=1}^cu_i=2t+1-p$,
simplifies to
\begin{equation}\label{eq:ineq_u2}
\begin{aligned}
2t(t-1)
&\le \sum_{i=1}^c u_i^2 + p^2 + 4p - 2cp - 4c(c-1) + (4c-5)(2t+1) \\
&\le \sum_{i=1}^c u_i^2 + p(p - 2c + 2) + (4c-5)(2t+1) \\
&\le \sum_{i=1}^c u_i^2 + 4c(c-1)^3 + (4c-5)(2t+1)
\end{aligned}
\end{equation}

Assume that there exists at most one line in $\Delta$ through $\infty$.
Then, by Lemma~\ref{lm:pre}(1),
for at most one $i^*\in \{1,2,\ldots,c\}$
we have $1+\kappa t\leq u_{i^*}\leq t$,
and $u_i < 1 + \kappa t$ for each $i\in\{1,2,\ldots,c\}\setminus \{i^*\}$.
Since $u_i\geq u_j$ implies that, for every $\epsilon>0$,
\[
(u_{i}+\epsilon)^2 + (u_j-\epsilon)^2 =
u_{i}^2+u_j^2+2\epsilon(u_{i}-u_j)+2\epsilon^2 >
u_{i}^2+u_j^2,
\]
one can see that
\begin{equation}\label{eq:sum_u}
\sum_{i=1}^c u_i^2 < t^2+(1+\kappa t)^2+(2t-t-\kappa t)^2.
\end{equation}

Combining Equations \eqref{eq:ineq_u2} and \eqref{eq:sum_u},
we obtain
\[
2t(t-1) < t^2+(1 + \kappa t)^2+(2t-t-\kappa t)^2 + 4c(c-1)^3 + (4c-5)(2t+1),
\]
which is violated if $t \geq 292$ for all $1 \le c \le 9$.

Thus, if $t\geq 292$ holds,
we obtain that the vertex $\infty$ lies in at least two lines.
Let $L_1$ and $L_2$ be two such lines.
Since a pair of nonadjacent vertices of $\Delta$
has at most $5$ common neighbours,
$|(L_1\cap L_2)\setminus \{\infty\}|\leq 4$ holds,
and hence $\Delta(\infty)\setminus (L_1\cup L_2)$ contains at most
\[
2t+1-(2\kappa t-4)=2(1-\kappa)t+5<\kappa t+1
\]
vertices.
This implies that $\infty$ lies in at most two lines.
Since $\infty$ was arbitrarily chosen in $\Delta$,
this shows the lemma for every vertex of $\Delta$.
\end{proof}

\begin{lemma}\label{lm:linesize}
Every line of $\Delta$ has size $t+1$.
\end{lemma}
\begin{proof}
Fix a vertex $\infty$ of $\Delta$ and,
following Lemma~\ref{lm:2lines},
let $L_1,L_2$ be the two lines through $\infty$.
Put $\ell_i := t+1-|L_i|$ ($i\in \{1,2\}$),
$R:=\Delta(\infty)\setminus (L_1\cup L_2)$, $r:=|R|$,
and $M:=\Delta(\infty)\cap L_1\cap L_2$, $m:=|M|$.
Without loss of generality, we assume that $\ell_1\leq \ell_2$.
Note that since each pair of nonadjacent vertices of $\Delta$ has at most $5$ common neighbours,
$m\leq 4$ follows.

It follows from Lemma~\ref{lm:pre}(3)
that the number of edges in $\Delta(\infty)$ equals $e_{\infty}:=t(t-1)$.
Let us estimate the number of edges in $\Delta(\infty)$
in terms of $\ell_1$, $\ell_2$, $m$, and $r$:
\begin{itemize}
\item \underline{Edges in $L_1\cup L_2$}.
As $L_i\setminus\{\infty\}$ is a clique,
it contains ${t-\ell_i\choose 2}$ edges;
then $L_1\cup L_2$ contains
${t-\ell_1\choose 2}+{t-\ell_2\choose 2}-{m\choose 2}$ edges.
\item \underline{Edges in $R$}.
The number of edges in $R$ is at most ${r\choose 2}$.
\item \underline{Edges between $R$ and $L_1\cup L_2$}.
Since each vertex in $R$
is adjacent to at most $4$ vertices in $L_i\setminus (M\cup \{\infty\})$
($i \in \{1, 2\}$),
the number of edges between $R$ and $L_1\cup L_2\setminus\{\infty\}$
is at most $2\cdot 4\cdot r$.
\item \underline{Edges between $L_1$ and $L_2$}.
Put $L_i':=L_i\setminus (M\cup \{\infty\})$.
We claim that the number $e_{12}$ of edges between $L_1'$ and $L_2'$
is at most $2(t-\ell_1)$.
Indeed,
let $s_i$ be the maximum number of edges
between a vertex $x\in L_i'$ and $L_j'$ ($\{i, j\} = \{1, 2\}$).
Then clearly
\begin{equation}\label{eq:e12-1}
e_{12}\leq \min\{s_1(t-\ell_1-m), s_2(t-\ell_2-m)\}.
\end{equation}
On the other hand,
no vertex $y\in L_j'\setminus \Delta(x)$
can be adjacent to more than $4-s_i-m$ vertices of $L_i'$,
as otherwise $|\Delta(x)\cap \Delta(y)|>s_i+1+m+(4-s_i-m)=5$ (recall that $\{\infty\} \cup M\subseteq \Delta(x)\cap \Delta(y)$),
a contradiction.
A vertex $y'\in L_j'\cap \Delta(x)$
can be adjacent to no more than $s_j-1$ vertices of $L_i'\setminus \{x\}$.
Hence,
\begin{align}
\nonumber e_{12} &\leq s_1 s_2 +
\min\{(t-\ell_2-m-s_1)(4-s_1-m), (t-\ell_1-m-s_2)(4-s_2-m)\}\\
& \leq s_1 s_2 +
\min\{(t-\ell_2-m)(4-s_1-m), (t-\ell_1-m)(4-s_2-m)\}\label{eq:e12-2}
\end{align}

Note that $t - \ell_1 - m \ge t - \ell_2 - m \ge \kappa t - 4 > 192 > 16 \ge s_1 s_2$,
so comparing Equations \eqref{eq:e12-1} and \eqref{eq:e12-2},
one can see that $e_{12}\leq 2(t-\ell_1)$.
\end{itemize}

Since $|\Delta(\infty)|=2t+1=(t-\ell_1)+(t-\ell_2)-m+r$ holds, one has $r=\ell_1+\ell_2+m+1$.

\begin{claim}\label{cl:nonzero}
Either $\ell_1=\ell_2=0$, or $\ell_1\geq 1$, $\ell_2\geq 1$.
\end{claim}
\begin{proof}
Without loss of generality, let $\ell_1=0$ and $\ell_2>0$. Then the clique $L_1$ has size $t+1$. By Lemma~\ref{lm:pre}(1),
there are no edges between $L_1\setminus\{\infty\}$ and $\Delta(\infty)\setminus L_1$; in particular, $m=0$ and $r=\ell_2+1$.
Thus, the number of edges in $\Delta(\infty)$ cannot exceed
\[
{t\choose 2}+{t-\ell_2\choose 2} + {\ell_2+1 \choose 2} + 4(\ell_2+1)=\ell_2^2 - \ell_2t + 5\ell_2 + t^2 - t + 4.
\]
which thus must be at least $e_{\infty}$. This implies $\ell_2^2 - \ell_2t + 5\ell_2 + 4\geq 0$, i.e.,
$t\leq\frac{\ell_2^2 + 5\ell_2 + 4}{\ell_2}=\ell_2 + 5 + \frac{4}{\ell_2}$,
a contradiction as $\ell_2=t+1-|L_2|$.
\end{proof}

Thus, assuming $\ell_1>0$ and $\ell_2>0$, the number of edges in $\Delta(\infty)$ does not exceed
\[
{t-\ell_1\choose 2}+{t-\ell_2\choose 2}-{m\choose 2} + {\ell_1+\ell_2+m+1 \choose 2} + 8(\ell_1+\ell_2+m+1) + 2(t-\ell_1),
\]
which hence must be at least $e_{\infty}$. This gives the inequality
\[
t(-\ell_1 - \ell_2 + 2)\geq -\ell_1^2 - \ell_1(\ell_2+m+7) - \ell_2^2 - \ell_2(m + 9) - 9m - 8,
\]
whence either $\ell_1+\ell_2>2$ and
\begin{equation}\label{eq:t}
t \leq {\ell_1^2 + \ell_1(\ell_2+m+7) + \ell_2^2 + \ell_2(m + 9) + 9m + 8
\over \ell_1 + \ell_2 - 2} =: d(\ell_1,\ell_2,m)
\end{equation}
or $\ell_1+\ell_2\leq 2$.

Suppose first that $\ell_1+\ell_2>2$.
It follows from Equation \eqref{eq:t} and $|L_i|\geq 1+\kappa t$ that
\begin{equation}\label{eq:l_over_t}
\frac{\ell_i}{d(\ell_1,\ell_2,m)}\leq \frac{\ell_i}{t} =
\frac{t+1-|L_i|}{t}\leq \frac{t+1-(1+\kappa t)}{t}=1-\kappa.
\quad (i \in \{1, 2\})
\end{equation}

By Claim~\ref{cl:nonzero}, one can see that
\[
\ell_1\geq 1, \ell_2\geq 1, \ell_1+\ell_2\geq 3, \ell_1,\ell_2\in \mathbb{N}, \text{~and~} m\in \{0,1,2,3,4\}.
\]
Subject to these constraints,
the inequalities defined by Equation \eqref{eq:l_over_t}
have a finite set $S$ of solutions in $(\ell_1,\ell_2,m)$,
which in turn shows that $\max_{S} d(\ell_1,\ell_2,m) = 88$.
By Equation \eqref{eq:t}, $t\leq 88$, a contradiction.

Thus, $\ell_1+\ell_2\leq 2$ holds.
It now follows from Claim~\ref{cl:nonzero}
that either $\ell_1=\ell_2=0$, or $\ell_1=\ell_2=1$.
Suppose the latter equality holds.
Then the cliques $L_1, L_2$ have size $t$.
Since $\Delta$ is connected and the vertex $\infty$ was chosen arbitrarily,
it follows from the above argument that every line in $\Delta$ has size $t$.
Since each vertex belongs to exactly two lines,
the number of lines in $\Delta$ equals $\frac{2|\Delta|}{t}=\frac{2(t^2+3t+2)}{t}\notin \mathbb{N}$,
a contradiction.
This shows the result.
\end{proof}

Proposition~\ref{pr:delta} immediately follows from Lemmas~\ref{lm:2lines}, \ref{lm:linesize}, the fact that every vertex in $\Delta$ has degree $2t+1$, and Lemma~\ref{lm:localTheta}.

This allows us to constrain the local graphs of $\G$
if the diameter $d$ is large enough.

\begin{proposition} \label{prop:locclas}
Let $\G$ be a distance-regular graphs
with classical parameters $(d, 2, 1, 2^d - 2)$ and diameter $d \ge 9$.
Then the local graphs of $\G$ are isomorphic to $\Lambda_{2^d-1}$.
\end{proposition}

\begin{proof}
By Corollary~\ref{cor:unique},
all local graphs of $\G$ have the same spectrum \eqref{eq:classpectrum},
which matches the spectrum \eqref{eq:spectrum} with $t = 2^d - 3 \ge 509$.
Let $\Delta = \G(u)$ be the local graph at some vertex $u$ of $\G$
-- this graph has $k(\G) = (2^d - 1)(2^d - 2) = (t+1)(t+2)$ vertices.
Since $c_2(\G) = 6$,
each pair of nonadjacent vertices of $\Delta$
has at most $5$ common neighbours in $\Delta$.
Then, by Proposition~\ref{pr:delta}, $\Delta$ is  isomorphic to $\Lambda_{t+2} = \Lambda_{2^d-1}$.
\end{proof}

\section{Locally $\Lambda_r$ graphs} \label{sec:loclm}

Recall that for a pair of vertices $u, v$ at distance $2$ in a graph $\G$,
the {\em $\mu$-graph} $M_\G(u, v)$ is the subgraph of $\G$
induced on the common neighbours of $u$ and $v$.
Furthemore,
an induced subgraph $\Delta$ of $\G$ is said to be {\em geodetically closed}
if for each pair of vertices $u, v$ of $\Delta$,
the graph $\Delta$ also contains every vertex $w$
lying on a shortest path between $u$ and $v$ in $\G$.
For a subset $U$ of the vertex set of $\G$,
the {\em geodetic closure} $G_\G(U)$
is the minimal geodetically closed subgraph of $\G$
that contains every vertex of $U$.

Let $\G$ be a locally $\Lambda_r$ graph with $c_2(\G) = 6$,
where $r \ge 4$.
We will investigate the $\mu$-graphs
and the geodetic closures of pairs of vertices at distance $2$ in $\G$.
To this end,
we will first introduce some definitions and basic results.

We note that the assumption on the local graphs of $\G$
implies that $k(\G) = r(r-1)$ and $a_1(\G) = 2r-3$,
and therefore $b_1(\G) = (r-1)(r-2)$
and $k_2(\G) = {r(r-1)^2 (r-2) \over 6}$.
Since the local graphs of $\G$ are isomorphic to $\Lambda_r$,
we will assume the existence of isomorphisms $\phi_u : \Lambda_r \to \G(u)$
for each vertex $u$ of $\G$.
Furthermore,
$\G$ has maximal cliques of sizes $3$ and $r$,
and we will call them {\em triangles} and {\em lines} in $\G$, respectively.
We also note that a maximal clique of $\G$ is uniquely determined
by each triple of its vertices,
and each pair of adjacent vertices of $\G$
lies in precisely one triangle and two lines.
In fact,
the lines and triangles containing a vertex $u$ of $\G$
are precisely
$\{u, (i, k)^{\phi_u} \mid 1 \le k \le r, k \ne i\}$ ($1 \le i \le r$),
$\{u, (k, j)^{\phi_u} \mid 1 \le k \le r, k \ne j\}$ ($1 \le j \le r$)
and $\{u, (i, j)^{\phi_u}, (j, i)^{\phi_u}\}$ ($1 \le i < j \le r$).

Consider a pair of vertices $(i, j), (k, \ell)$
($1 \le i, j, k, \ell \le r$, $i \ne j$, $k \ne \ell$)
of $\Lambda_r$ at distance $2$,
i.e., $i \ne k$, $j \ne \ell$, and $(i, j) \ne (\ell, k)$.
Note that the vertices $(i, j)$ and $(k, \ell)$
have three common neighbours in $\Lambda_r$ if $i = \ell$ or $j = k$,
and two common neighbours otherwise.
Let $S = \{u, v, w\}$ be an independent set of $\Lambda_r$
such that there exist distinct vertices $x, y$
that are adjacent to all vertices in $S$.
We will call such a set $S$ a {\em semibiclique} of $\Lambda_r$.
By the previous observation,
we see that we have $S = \{(i, j), (j, k), (k, i)\}$
for some distinct $1 \le i, j, k \le r$,
and the common neighbours of each pair of vertices of $S$
form another semibiclique $S' = \{(j, i), (k, j), (i, k)\}$,
called the {\em complementary semibiclique of $S$},
such that the subgraph of $\Lambda_r$ induced on $S \cup S'$
is isomorphic to $K_{3, 3}$.
Note also that $S'' = S$.
From this we also see that a pair of vertices of $\Lambda_r$
lies in at most one semibiclique.

Let us now give some lemmas which will be useful later.

\begin{lemma} \label{lem:two}
Let $u$ be a vertex of $\G$, and $v, w$ nonadjacent vertices of $\G(u)$.
Then precisely two common neighbours of $v$ and $w$ in $\G(u)$
lie in lines containing $u$ and $v$.
\end{lemma}

\begin{proof}
If $v$ and $w$ have three common neighbours in $\G(u)$,
then two of them lie in lines containing $u$ and $v$
and one lies in the triangle containing $u$ and $v$.

Otherwise,
we may write $v = (i, j)^{\phi_u}$ and $w = (k, \ell)^{\phi_u}$
for some distinct $1 \le i, j, k, \ell \le r$.
Then,
the common neighbours of $v$ and $w$ in $\G(u)$
are $(i, \ell)^{\phi_u}$ and $(k, j)^{\phi_u}$,
and both $\{u, v, (i, \ell)^{\phi_u}\}$ and
$\{u, v, (k, j)^{\phi_u}\}$ determine lines.
\end{proof}

\begin{lemma} \label{lem:three}
Let $\{u, v, w\}$ be a triangle in $\G$.
Then a vertex $x$ of $\G(u)$
has three common neighbours with $v$ in $\G(u)$
precisely when it lies on a line containing $u$ and $w$.
\end{lemma}

\begin{proof}
We note that $v = (i, j)^{\phi_u}$ and $w = (j, i)^{\phi_u}$
for some distinct $1 \le i, j \le r$.
Then $x$ has three common neighbours with $v$ in $\G(u)$
precisely when $x = (k, i)^{\phi_u}$ or $x = (j, k)^{\phi_u}$
for some $1 \le k \le r$ with $k \ne i, j$,
i.e., when there is an $(r-1)$-clique in $\G(u)$ containing $w$ and $x$,
or equivalently, $\{u, w, x\}$ determines a line in $\G$.
\end{proof}

\begin{lemma} \label{lem:bij}
Let $u$ be a vertex of $\G$ and $S$ a semibiclique in $\G(u)$
with corresponding complementary semibiclique $S'$ in $\G(u)$.
Then there exists a bijective mapping $\psi: S \to S'$
such that for all $v \in S$,
the set $\{u, v, v^\psi\}$ is a triangle.
\end{lemma}

\begin{proof}
We have $S = \{(i, j)^{\phi_u}, (j, k)^{\phi_u}, (k, i)^{\phi_u}\}$
for some distinct $1 \le i, j, k \le r$,
and then also $S' = \{(j, i)^{\phi_u}, (k, j)^{\phi_u}, (i, k)^{\phi_u}\}$.
Therefore,
the map
$\psi: (a, b)^{\phi_u} \mapsto (b, a)^{\phi_u}$ ($(a, b)^{\phi_u} \in S$)
gives the desired bijective mapping.
\end{proof}

\begin{lemma} \label{lem:semibc}
Let $u$ and $v$ be vertices of $\G$ at distance $2$
contained in a semibiclique $S$ in $\G(x)$,
where $x$ is a common neighbour of $u$ and $v$ in $\G$.
Then the complementary semibiclique $S'$ of $S$ in $\G(x)$
is also a semibiclique in $\G(u)$.
\end{lemma}

\begin{proof}
By assumption,
we can write $x = (i, j)^{\phi_u}$ for some distinct $1 \le i, j \le r$.
Then the set $S'$ contains a vertex
from each of the two lines and the triangle containing $u$ and $x$,
so we have $S' = \{(j, i)^{\phi_u}, (k, j)^{\phi_u}, (i, k')^{\phi_u}\}$
for some $1 \le k, k' \le r$ ($k, k' \ne i, j$).

By Lemma~\ref{lem:bij} and since $\{u, x, (j, i)^{\phi_u}\}$ is a triangle,
we may assume without loss of generality
that $\{v, x, (k, j)^{\phi_u}\}$ is a triangle.
Since $\{v, x, (j, i)^{\phi_u}\}$ determines a line,
by Lemma~\ref{lem:three},
the vertices $(k, j)^{\phi_u}$ and $(j, i)^{\phi_u}$
have three common neighbours in $\G(v)$.
As $c_2(\G) = 6$,
it follows that these common neighbours
are precisely the ones in $\G(u)$,
namely, $x = (i, j)^{\phi_u}$, $(j, k)^{\phi_u}$ and $(k, i)^{\phi_u}$.
Similarly,
the common neighbours of $(i, k')^{\phi_u}$ and $(j, i)^{\phi_u}$ in $\G(v)$
are $x$, $(j, k')^{\phi_u}$ and $(k', i)^{\phi_u}$.
Of the two lines containing $v$ and $(j, i)^{\phi_u}$,
one thus contains $x$,
while, by Lemma~\ref{lem:two},
the other contains one of the common neighbours in $\G(v)$
of $(j, i)^{\phi_u}$ with each of $(k, j)^{\phi_u}$ and $(i, k')^{\phi_u}$.
If $k \ne k'$,
then the latter line contains
either $(j, k)^{\phi_u}$ and $(j, k')^{\phi_u}$
or $(k, i)^{\phi_u}$ and $(k', i)^{\phi_u}$,
and must thus also contain $u$, contradiction.
Therefore, $k = k'$ holds, and $S'$ is a semibiclique in $\G(u)$.
\end{proof}

We will now identify the $\mu$-graphs in $\G$
along with the corresponding geodetic closures.

\begin{proposition} \label{prop:three}
Let $u$ be a vertex of $\G$.
If $v$ is a vertex of $\G_2(u)$
such that $u$ and $v$ have three common neighbours in $\G(x)$,
where $x$ is a common neighbour of $u$ and $v$ in $\G$,
then $M_\G(u, v) \cong K_{3,3}$.
Moreover,
there exist precisely ${r(r-1)(r-2) \over 3}$ such vertices $v$ in $\G_2(u)$,
and then the geodetic closure $G_\G(\{u, v\})$ is isomorphic to $K_{3,3,3}$
and contains two vertices $v, w$ from $\G_2(u)$,
for which $M_\G(u, v) = M_\G(u, w)$ and $G_\G(\{u, v\}) = G_\G(\{u, w\})$ hold.
\end{proposition}

\begin{proof}
The vertex $v$ has three common neighbours with $u$ in $\G(x)$
precisely when there exists a semibiclique $S = \{u, v, w\}$ in $\G(x)$.
By Lemma~\ref{lem:semibc},
the complementary semibiclique $S'$ of $S$ in $\G(x)$
is also a semibiclique in each of $\G(u)$, $\G(v)$ and $\G(w)$
-- in particular,
we can write $x = (i, j)^{\phi_u}$
and $S' = \{(j, i)^{\phi_u}, (k, j)^{\phi_u}, (i, k)^{\phi_u}\}$
for some distinct $1 \le i, j, k \le r$.

Since $S'' = S$ is the complementary semibiclique of $S'$ in $\G(x)$,
Lemma~\ref{lem:semibc} implies that
$S$ is also a semibiclique in $\G((j, i)^{\phi_u})$,
and its complementary semibiclique $S^*$ in $\G((j, i)^{\phi_u})$
is also a semibiclique in each of $\G(u)$, $\G(v)$ and $\G(w)$.
The set of common neighbours of $u, v, w$ in $\Gamma$
is then the set $S' \cup S^*$ of size $c_2(\G) = 6$,
so $S^*$ is the complementary semibiclique of $S'$ in $\G(u)$.
We thus have $S^* = \{x, (j, k)^{\phi_u}, (k, i)^{\phi_u}\}$.
Therefore,
the $\mu$-graph $M_\G(u, v) = M_\G(v, w) = M_\G(w, u)$
has vertex set $S' \cup S^*$ and is isomorphic to $K_{3, 3}$.
Furthermore,
the geodetic closure $G_\G(\{u, v\}) = G_\G(\{v, w\}) = G_\G(\{w, u\})$
has vertex set $S \cup S' \cup S^*$ and is isomorphic to $K_{3, 3, 3}$;
in particular, it contains the vertices $v$ and $w$ from $\G_2(u)$.

Suppose now that for some $v \in \G_2(u)$,
we have $M_\G(u, v) \cong K_{3, 3}$.
Since this is a cubic graph,
it follows that $u$ and $v$ have three common neighbours in $\G(x)$
for each of the common neighbours $x$ of $u$ and $v$ in $\G$.
We note that for every vertex $x$ of $\G(u)$,
the vertex $v$ can be chosen
from one of the semibicliques of $\G(x)$ containing $u$.
The above argument thus shows that
for every complementary pair $S', S^*$ of semibicliques in $\G(u)$
there exist two vertices $v$ in $\G_2(u)$
such that the set of common neighbours of $u$ and $v$ is $S' \cup S^*$.
Since there is such a pair of semibicliques
for each choice of $1 \le i < j < k \le r$,
it follows that there are
$2 \cdot {r \choose 3} = {r(r-1)(r-2) \over 3}$ vertices $v$ in $\G_2(u)$
such that $M_\G(u, v) \cong K_{3,3}$.
\end{proof}

\begin{proposition} \label{prop:two}
Let $u$ be a vertex of $\G$.
If $v$ is a vertex of $\G_2(u)$
such that $u$ and $v$ have two common neighbours in $\G(x)$,
where $x$ is a common neighbour of $u$ and $v$ in $\G$,
then $M_\G(u, v) \cong C_6$.
Moreover,
there exist precisely ${r(r-1)(r-2)(r-3) \over 6}$
such vertices $v$ in $\G_2(u)$,
and then the geodetic closure $G_\G(\{u, v\})$ is isomorphic to $H_2(2, 2)$
and contains six vertices from $\G_2(u)$;
also,
for each vertex $z$ of $G_\G(\{u, v\})$ in $\G_2(u)$,
$G_\G(\{u, v\}) = G_\G(\{u, z\})$ holds.
\end{proposition}

\begin{proof}
Since $u$ and $v$ have two common neighbours in $\G(x)$,
we must have $M_\G(u, v) \not\cong K_{3,3}$,
and by Proposition~\ref{prop:three},
it follows that $u$ and $v$ have two common neighbours in $\G(y)$
for every common neighbour $y$ of $u$ and $v$,
so $M_\G(u, v)$ is a union of cycles of total length $c_2(\G) = 6$.
As we can write $x = (h, \ell)^{\phi_u}$
for some distinct $1 \le h, \ell \le r$,
Lemma~\ref{lem:two} implies that
the two common neighbours of $u$ and $v$ in $\G(x)$
are $(h, m)^{\phi_u}$ and $(j, \ell)^{\phi_u}$
for some distinct $1 \le j, m \le r$ with $j, m \ne h, \ell$.
Since the latter two vertices are nonadjacent,
it follows that $M_\G(u, v)$ is triangle-free and is therefore the $6$-cycle
$((h, \ell)^{\phi_u}, (h, m)^{\phi_u}, (i, m)^{\phi_u}, (i, k)^{\phi_u},
(j, k)^{\phi_u}, (j, \ell)^{\phi_u})$
for some distinct $1 \le i, k \le r$ with $i, k \ne h, j, \ell, m$.
Proposition~\ref{prop:three} also implies that there are precisely
$k_2(\G) - {r(r-1)(r-2) \over 3} = {r(r-1)(r-2)(r-3) \over 6}$
vertices $v$ such that $M_\G(u, v) \cong C_6$.

Note that each pair of consecutive vertices in the above $6$-cycle
determines a line together with $u$,
and the three pairs of such lines
corresponding to a pair of antipodal edges in the $6$-cycle
have a second common vertex,
namely $(h, k)^{\phi_u}$, $(i, \ell)^{\phi_u}$ and $(j, m)^{\phi_u}$.
Similarly,
each pair of consecutive vertices in the above $6$-cycle
also determines a line together with $v$,
and we denote the common vertices other than $v$
of the lines determined by $\{v, (h, \ell)^{\phi_u}, (h, m)^{\phi_u}\}$
and $\{v, (i, k)^{\phi_u}, (j, k)^{\phi_u}\}$,
by $\{v, (i, k)^{\phi_u}, (i, m)^{\phi_u}\}$
and $\{v, (h, \ell)^{\phi_u}, (j, \ell)^{\phi_u}\}$,
and by $\{v, (j, k)^{\phi_u},\allowbreak (j, \ell)^{\phi_u}\}$
and $\{v, (h, m)^{\phi_u}, (i, m)^{\phi_u}\}$,
by $w$, $w'$ and $w''$,
respectively.

We will now construct the geodetic closure $G_\G(\{u, v\})$.
Proposition~\ref{prop:three} and the above argument imply that
for each pair of nonadjacent vertices $y, z$ of $G_\G(\{u, v\})$,
we have $M_\G(y, z) \cong C_6$.
Note that since the set
$\{(h, k)^{\phi_u}, (h, \ell)^{\phi_u}, (h, m)^{\phi_u}\}$
determines a line containing $u$,
it follows that $w$ is not adjacent to $(h, k)^{\phi_u}$.
Since $w$ is adjacent to $v$,
it is therefore not adjacent to $(i, m)^{\phi_u}$ and $(j, \ell)^{\phi_u}$
and is then adjacent to $(i, \ell)^{\phi_u}$ and $(j, m)^{\phi_u}$.
Similarly, $w'$ is adjacent to $(h, k)^{\phi_u}$ and $(j, m)^{\phi_u}$,
and $w''$ is adjacent to $(h, k)^{\phi_u}$ and $(i, \ell)^{\phi_u}$.
Now, the common neighbours of $(h, k)^{\phi_u}$ and $w$,
besides $(h, \ell)^{\phi_u}, (h, m)^{\phi_u}, (i, k)^{\phi_u}, (j, k)^{\phi_u}$,
are two vertices $v'$ and $v''$
such that $v'$ is adjacent to $(h, \ell)^{\phi_u}$ and $(j, k)^{\phi_u}$,
and $v''$ is adjacent to $(h, m)^{\phi_u}$ and $(i, k)^{\phi_u}$.
Since the set $\{(h, k)^{\phi_u}, (h, \ell)^{\phi_u}, (h, m)^{\phi_u}\}$
determines a line containing $u$,
which is not adjacent to $v'$ or $w'$,
it follows that the sets $\{(h, k)^{\phi_u}, (h, \ell)^{\phi_u}, v'\}$
and $\{(h, k)^{\phi_u}, (h, \ell)^{\phi_u}, w'\}$
determine the same line.
This and similar arguments show that both $v'$ and $v''$
are adjacent to each of $(i, \ell)^{\phi_u}, (j, m)^{\phi_u}, w', w''$.
Finally,
by examining the common neighbours of $(i, \ell)^{\phi_u}$ and $w'$,
it follows that $v'$ is also adjacent to $(i, m)^{\phi_u}$
and $v''$ is also adjacent to $(j, \ell)^{\phi_u}$.
Thus, the geodetic closure of $\{u, v\}$ has vertex set
\begin{align*}
\{&{} u, v, v', v'', w, w', w'', (h, k)^{\phi_u}, (h, \ell)^{\phi_u},
(h, m)^{\phi_u}, \\ &{} (i, k)^{\phi_u}, (i, \ell)^{\phi_u}, (i, m)^{\phi_u},
(j, k)^{\phi_u}, (j, \ell)^{\phi_u}, (j, m)^{\phi_u}\}
\end{align*}
and is isomorphic to $H_2(2, 2) \cong \overline{K_4 \square K_4}$;
in particular, it contains six vertices from $\G_2(u)$.
Since for each such vertex $z$,
the graph $M_\G(u, z)$ is a subgraph of $G_\G({u, v})$ isomorphic to $C_6$,
it follows that $G_\G({u, v}) = G_\G({u, z})$.
\end{proof}

We may now derive some divisibility conditions,
which will allow us to rule out an infinite family
of distance-regular graphs with classical parameters.

\begin{corollary} \label{cor:div}
Let $\G$ be a locally $\Lambda_r$ graph with $c_2(\G) = 6$,
where $r \ge 4$,
and let $n$ denote the number of vertices of $\G$.
Then, the following hold.
\begin{enumerate}[(a)]
\item $r \not\equiv 4, 5, 7, 8 \pmod{9}$.
\item $n \cdot {r \choose 3}$ is divisible by $9$.
\item $n \cdot {r \choose 4}$ is divisible by $8$.
\end{enumerate}
\end{corollary}

\begin{proof}
Let $u$ be a vertex of $\G$.
By Proposition~\ref{prop:two},
there are precisely ${r(r-1)(r-2)(r-3) \over 6}$ vertices $v$ of $\G_2(u)$
such that $M_\G(u, v) \cong C_6$,
and the relation on $\G_2(u)$ of being in the same geodetic closure with $u$
is an equivalence relation with equivalence classes of size $6$.
It then follows that ${r(r-1)(r-2)(r-3) \over 36}$ is an integer,
which is true whenever $r \not\equiv 4, 5, 7, 8 \pmod{9}$.
This proves part (a).

Let us now count the number of geodetically closed $K_{3,3,3}$ subgraphs in $\G$.
By Proposition~\ref{prop:three},
each vertex of $\G$ lies on
${r(r-1)(r-2) \over 6} = {r \choose 3}$ such subgraphs,
and each such subgraph contains $9$ vertices.
Therefore,
there are precisely ${n \cdot {r \choose 3} \over 9}$
geodetically closed $K_{3,3,3}$ subgraphs in $\G$,
from which part (b) follows.

Similarly,
we now count the number of geodetically closed $H_2(2, 2)$ subgraphs in $\G$.
By the argument above,
each vertex of $\G$ lies on
${r(r-1)(r-2)(r-3) \over 36} = {2 \over 3} {r \choose 4}$ such subgraphs,
and each such subgraph contains $16$ vertices.
Therefore,
there are precisely ${n \cdot {r \choose 4} \over 3 \cdot 8}$
geodetically closed $H_2(2, 2)$ subgraphs in $\G$,
which proves part (c).
\end{proof}

\begin{proposition} \label{prop:loc}
A distance-regular graph
with classical parameters $(d, 2, 1, 2^d-2)$ and $d \ge 3$
is not a locally $\Lambda_{2^d-1}$ graph.
\end{proposition}

\begin{proof}
Let $\G$ be a distance-regular graph
with classical parameters $(d, 2, 1, 2^d-2)$,
where $d \ge 3$.
By~\cite[Corollary~8.4.4]{BCN}, $\G$ has $n = (2^d - 1)^d$ vertices.
Thus,
$$
n \cdot {2^d - 1 \choose 4} = {n(2^d - 1)(2^d - 2)(2^d - 3)(2^d - 4) \over 24} =
{n(2^d - 1)(2^{d-1} - 1)(2^d - 3)(2^{d-2} - 1) \over 3}
$$
is an odd integer and therefore not divisible by $8$.
Since $c_2 = 6$,
Corollary~\ref{cor:div}(c) then implies
that $\G$ is not a locally $\Lambda_{2^d-1}$ graph.
\end{proof}

\begin{theorem}
A distance-regular graph
with classical parameters $(d, 2, 1, 2^d-2)$ and $d \ge 9$
does not exist.
\end{theorem}

\begin{proof}
Let $\G$ be a distance-regular graph
with classical parameters $(d, 2, 1, 2^d-2)$ and diameter $d \ge 9$.
By Proposition~\ref{prop:locclas},
the local graphs of $\G$ are isomorphic to $\Lambda_{2^d-1}$,
which contradicts Proposition~\ref{prop:loc}.
\end{proof}

\section*{AI use statement}
During the preparation of this work the authors used AI
in order to simplify the proof of Theorem~\ref{thm:nonext}\footnote{%
The original proof reduced the problem of ruling out the cases $u = 3, 4, 5$ to a Pell-like equation and then used computer algebra to solve it.
}.
The authors reviewed and edited its output as needed.
All other mathematical content of the paper is solely due to the authors.

\section*{Acknowledgements}
Janoš Vidali is supported in part
by the Slovenian Research and Innovation Agency
(research program P1-0285 and research projects
J1-1691, J1-1692, J1-2452, J1-3002, N1-0216 and N1-0218).

\end{document}